\documentclass{templt}

\newtheorem{thm}{Theorem} 

\newtheorem{cor}[thm]{Corollary}
\newtheorem{rem}[thm]{Remark}
\newtheorem{defn}[thm]{Definition}
\newtheorem{prop}[thm]{Proposition}
\newtheorem{ex}[thm]{Example}

\usepackage{hyperref}
\usepackage{orcidlink}
\numberwithin{equation}{section}

\def\C{\mathbb{C}}
\begin{document}

\title[Operators On Smooth Banach Spaces]{Normal and Unitary Operators on Smooth Banach Spaces via the Duality Mapping}

\author[Mohammed Shameem]{Mohammed Shameem \orcidlink{0009-0003-2148-3750}}
\address{Department of Mathematics, National Institute of Technology Calicut, Kozhikode, Kerala, India 673601
ORCID: \href{https://orcid.org/0009-0003-2148-3750}{0009-0003-2148-3750}}
\curraddr{}
\email{shameemmaths@gmail.com}
\thanks{}

\author[Deepesh K. P.]{Deepesh K. P. \orcidlink{0009-0002-4790-9425}}
\address{Department of Mathematics, National Institute of Technology Calicut, Kozhikode, Kerala, India 673601
ORCID: \href{https://orcid.org/0009-0002-4790-9425}{0009-0002-4790-9425}}
\curraddr{}
\email{deepeshkp@nitc.ac.in}
\thanks{}


\date{}

\dedicatory{}

%    "Communicated by" -- provide editor's name; required.
% \commby{}

\begin{abstract}
This article extends the classical notions of normal and unitary operators from Hilbert spaces to smooth Banach spaces through the duality mapping. The class of normal operators introduced in this setting contains, as a special case, the self-adjoint operators on smooth Banach spaces defined via the duality mapping. We establish several characterisations of these operators and investigate their fundamental properties. In particular, we examine various spectral properties of normal operators on smooth Banach spaces. These results demonstrate that many fundamental aspects of the classical theory of normal and unitary operators on Hilbert spaces persist in the smooth Banach space setting.
\end{abstract}

\maketitle

\section{Introduction}\label{intro}
The study of special classes of operators like the self-adjoint, normal and unitary operators on Hilbert spaces, particularly their spectral properties, has a long and distinguished history and continues to be an active area of research \cite{carvajal2012operators, limaye1996functional}. Motivated by these operators, several authors have sought to extend these notions beyond Hilbert spaces, particularly to the setting of Banach spaces \cite{shtraus1978,wojcik2013}. In this direction, Stampli (1969) introduced the concepts of adjoint-abelian and iso-abelian operators on Banach spaces, which are analogous to self-adjoint and unitary operators on Hilbert spaces, respectively \cite{stampfli1969}. Later in $2020$, García-Pacheco introduced a notion of self-adjoint operators on Banach spaces using the duality mapping \cite{banachSA}. %Likewise, the concept of numerical range has been generalized to Banach spaces through continuous linear functionals \cite{duncan1971}.
\vskip0.25cm
Motivated by these developments, we introduce a new notion of normal operators on smooth Banach spaces, which is a natural extension of the classical notion of normal operators on Hilbert spaces. This class includes all self-adjoint operators on smooth Banach spaces defined in \cite{banachSA}. Our defintion naturally induces the corresponding notion of unitary operators, extending the classical concept to this framework. Several characterizations of these classes of operators are established in this article.
\vskip0.25cm
We further investigate the spectral properties of normal operators on smooth Banach spaces. In particular, we establish that the spectrum and the approximate spectrum coincide for these normal operators.  We prove that the eigenvectors corresponding to distinct eigenvalues of a normal operator are Birkhoff--James orthogonal, providing an analogue of the corresponding Hilbert space result. We also show that the eigenspectrum of a normal operator is countable. We show that the spectrum of a unitary operator is contained in the unit circle and derive several other spectral properties. These properties establishes the naturality of the extensions.
\vskip0.25cm

 The paper is organized as follows. In Section 2, we introduce the notations and preliminary results required for a smooth reading of the paper. Section 3 is devoted to the introduction and characterization of normal and unitary operators on smooth Banach spaces. In Section 4, we investigate the spectral properties of these operators, and we conclude the article with a few related open problems.

\section{Notations and Preliminaries}\label{prelim}
Throughout this article, $X$ and $Y$ denote Banach spaces over the complex field $\mathbb{C}$, and $B(X,Y)$ denotes the space of all bounded linear operators from $X$ into $Y$. When $X=Y$, we write $B(X)$ for $B(X,X)$, while $X^*:=B(X,\mathbb{C})$ denotes the space of all bounded linear functionals on $X$. The unit sphere of $X$, given by $\{x\in X:\|x\|=1\}$, is denoted by $S_X$, and $H$ denotes a complex Hilbert space. Here, the term operator refers to a bounded linear operator, and $\mathbb{N}$ denotes the set of all natural numbers.
\vskip0.25cm

For an operator $T\in B(X)$, the spectrum $\sigma(T)$ is defined as the set of all $\lambda\in\mathbb{C}$ for which $T-\lambda I$ does not have a bounded inverse. The eigenspectrum $\sigma_{eig}(T)$ consists of those $\lambda\in\mathbb{C}$ for which $T-\lambda I$ is not injective, while the approximate spectrum $\sigma_{app}(T)$ consists of those $\lambda\in\mathbb{C}$ for which $T-\lambda I$ is not bounded below. The spectral radius of $T$, denoted by $\rho(T)$, is given by
$
\rho(T)=\sup\{|\lambda|:\lambda\in\sigma(T)\}.
$
\vskip0.25cm
Two important quantities associated with $T\in B(X)$ are the operator norm $\|T\|$ and the minimum modulus $\mu(T)$, defined respectively by

$$
\|T\|=\sup\{\|Tx\|:x\in S_X\} \,\, \mbox{and}\,\, \mu(T)=\inf\{\|Tx\|:x\in S_X\}.
$$
\vskip0.25cm
Another important set, closely connected to the spectrum of an operator $T \in B(X)$ is the numerical range $W(T)$, which is given by
$$
W(T) = \{ x^*(T(x)) : x \in S_X, \, x^* \in S_{X^*},\, x^*(x) = 1 \}.
$$
Some quantities associated with the numerical range of an operator $T$ are the numerical radius $r(T)$ and the Crawford number $c(T)$ given by
$$
r(T) = \sup \{ |\lambda| : \lambda\in W(T) \}\,\, \mbox{and}\,\, c(T) = \inf \{ |\lambda| : \lambda\in W(T) \}.
$$
These quantities represent, respectively, the upper and lower bounds of the numerical range in modulus and are well studied in the literature \cite{crawford2024,crawfordattainset,crawfordattainsetc,sd2026}.\vskip0.25cm

The transpose of an operator $T \in B(X,Y)$ is defined as the linear operator $T^{\prime}$ from $Y^*$ to $X^*$ defined by 
$$
(T^{\prime}(y))(x) = y(Tx), \quad y\in Y^*,\, x \in X,
$$
and for an operator $T \in B(H_1,H_2)$, where $H_1,H_2$ are Hilbert spaces, the adjoint is defined as the operator $T^*$ from $H_2$ to $H_1$ defined by
$$
\langle Tx, y \rangle = \langle x, T^*y \rangle,\, x \in H_1,\, y\in H_2.
$$

A Banach space $X$ is said to be smooth if, for every $x\in S_X$, there exists a unique $x^*\in S_{X^*}$ such that $x^*(x)=1$ (\cite{dragomir2003}). Consequently, for a smooth Banach space $X$, one can define a duality map $J : X \to X^*$ satisfying the conditions
$$
\|J(x)\|=\|x\| \,\,\mbox{and}\,\, J(x)(x)=\|x\|^2\
$$
for every $x\in X$. It is well known (\cite{dragomir2003}) that, for a smooth Banach space $X$, the mapping $J$ is surjective if and only if $X$ is reflexive, and injective if and only if $X$ is strictly convex. Thus, for an operator $T$ on a reflexive strictly convex smooth Banach space, $J^{-1}T'J$ is well defined, and it is denoted by $T^*$ here, as it coincides with the adjoint operator $T^*$ in the Hilbert space settings.
\vskip0.25cm

Given two elements $x, y$ in a Banach space $X$, we say that $x$ is \emph{Birkhoff--James orthogonal} to $y$, denoted by $x \perp_B y$, if
$
\|x + \lambda y\| \geq \|x\|$ for all $\, \lambda \in \mathbb{C}
$ (\cite{mal}). In a smooth Banach space $X$, the Birkhoff--James orthogonality can be characterized in terms of the duality mapping $J$. That is, for $x,y\in X,$ $x\perp_B y$ if and only if $J(x)(y)=0$.
\vskip0.25cm
In the setting of smooth Banach spaces, a notion of self-adjoint operator was introduced in \cite{banachSA} that is analogous to the classical self-adjoint operator in the Hilbert space setting, as defined below.
\begin{defn}(\cite{banachSA})\label{SA}
    Let $X$ be a smooth Banach space. An operator $T \in B(X)$ is said to be \textit{self-adjoint} if $T' \circ J = J \circ T$.

\end{defn}
 
Several properties of these operators were studied and many results on them were established in the article \cite{banachSA}.

\section{Classes of normal operators and unitary operators}

On a Hilbert space, it is well known that the duality map $J$ is invertible and linear. Consequently, as mentioned before, for an operator $T$ on a Hilbert space, the operator $J^{-1}T'J$ coincides with the adjoint operator $T^*$ and hence, $T$ becomes self-adjoint if and only if $T = J^{-1}T'J$. A class of self adjoint operators on smooth Banach spaces was defined in \cite{banachSA} using this expression (as given in Definition \ref{SA}) and many properties analogous to the classical self-adjoint operators were established therein. Though the map $J$ is invertible when $X$ is a reflexive strictly convex smooth Banach space, due to the nonlinearity of the duality map, extending the notion of these special types of operators defined using the adjoint operator to the Banach space setting is a challenging problem, especially for normal operators (\cite{banachSA}). 
\vskip0.25cm

On a Hilbert space, it is well-known that a bounded linear operator $T$ is normal if and only if $\|Tx\|=\|T^*x\|$ for all $x \in X$ (\cite{limaye1996functional}). %, and it is unitary if $\|Tx\|=\|T^*x\|=\|x\|$ for all $x \in X$ [cite]. 
In view of this characterization, one can define a bounded linear operator $T$ on a reflexive strictly convex smooth Banach space $X$ to be \emph{normal} if $\|Tx\|=\|J^{-1}T'Jx\|$ for all $x \in X$. Further, $J$ being an isometry, this characterization can be extended to arbitrary smooth Banach spaces, where the invertibility of the duality map is no longer required.
\begin{defn}\label{normal-defn}
    Let $X$ be a smooth Banach space. Then $T \in B(X)$ is said to be 
    \begin{enumerate}
        \item normal if $\|Tx\| = \|T'Jx\| \, \mbox{ for all } x \in X.$
        \item unitary if $\|Tx\| = \|T'Jx\|=\|x\| \, \mbox{ for all } x \in X.$
    \end{enumerate}
    
\end{defn}

Note that in the case of a normal operator $T$ on a Hilbert space, since $T^* $ is equal to $J^{-1}T^{\prime}J$, the characterization $\|Tx\|=\|T^*x\|$ is equivalent to $\|Tx\| = \|T'Jx\|$ for all $x\in X$. Thus, the Definition \ref{normal-defn} naturally extends the existing concept of normal operator to smooth Banach space settings.

It is obvious that unitary operators are normal. When $T$ is a self-adjoint operator, due to the equality
\(\|T^{\prime}J(x)\|=\|JT(x)\|=\|Tx\|,\) 
it follows that every self-adjoint operator on $X$ is also normal. 
\vskip0.25cm
As in the Hilbert space setting, the scalar multiple of a special type operator need not be always of the same type.

\begin{prop}
 If $T\in B(X)$ is normal, then $\lambda T$ is normal for every $\lambda\in\mathbb{C}$. Moreover, if $T$ is unitary, then $\lambda T$ is unitary if and only if $|\lambda|=1$, while if $T$ is self-adjoint, then $\lambda T$ is self-adjoint if and only if $\lambda\in\mathbb{R}$.
\end{prop}
\begin{proof}
    For $x \in X$, $\lambda \in \C$ and $T$ normal, we get
    $$
    \|\lambda Tx\| = \|T(\lambda x)\|=\|T'J(\lambda x)\| = \|\overline{\lambda}T'Jx\|=\|(\lambda T)'Jx\|.
    $$ 
The unitary and self-adjoint cases are clear from their definitions.
\end{proof}

Though the sum of two normal operators need not be normal even if they commute (see Example \ref{normal_sum}), we have the following result on the transpose of a normal operator.

\begin{prop}\label{normal_transp}
    Let $X$ be reflexive smooth Banach space. Then $T \in B(X)$ is normal if and only if $T'$ is normal.
\end{prop}
\begin{proof}
    Let $T\in B(X)$ be a normal operator and $\tilde{J}$ be the duality map from $X^*$ to $X^{**}$. We need to show that $\|T'f\|=\|T''\tilde{J}f\|$, for all $f \in X^*$. Since $X$ is reflexive, for $f \in X^*$, there exists $x \in X$ such that $\|T'f\|=\|T'Jx\|=\|Tx\|$. Since $\tilde{J}Jx = \phi(x)$, where $\phi:X \rightarrow X^{**}$ is the canonical embedding map (\cite{banachSA}), we get $\|T''\tilde{J}f\|=\|T''\tilde{J}Jx\|=\|T''\phi(x)\|$. Now for $x\in X$ and $g \in X^*$, we have,
    $$
    (T''\phi(x))(g) = \phi(x)T'g = (T'g)(x) = \phi(Tx)(g). 
    $$
    Since $\phi$ is an isometry, $\|T''\tilde{J}f\| = \|\phi(Tx)\|=\|Tx\|$, and hence $T'$ is normal. 

When $X$ is reflexive, $T$ and $T''$ are equivalent. So if $T'$ is normal, applying the previous part, $T''$ is normal, which implies that $T$ is also normal.
\end{proof}

An similar result holds for unitary operators on reflexive smooth Banach spaces. 
\begin{prop}
    Let $X$ be a reflexive smooth Banach space and $T \in B(X)$. Then $T$ is unitary on $X$ if and only if $T'$ is a unitary operator on $X^*$.
\end{prop}
\begin{proof}
    Let $T\in B(X)$ be unitary operator. Then $T'$ is normal (Proposition \ref{normal_transp}), and hence,
    \[
        \|T'f\|=\|T''\tilde{J}f\|,
    \] for each $f \in X^*$, where $\tilde{J}$ is duality map on $X^*$. Since $X$ is reflexive, for any $f \in X^*$, there exists $x \in X$ such that $Jx=f$. Then
    \[
        \|T'f\|=\|T'Jx\|=\|x\|=\|Jx\|=\|f\|.    
    \] Thus $T'$ is unitary. Since $X$ is reflexive, the converse part is clear.
\end{proof}

Similar to the classical characterization of unitary operators in the Hilbert space setting, we characterize unitary operators on reflexive smooth Banach spaces.

\begin{thm}\label{unitary_charc}
    Let $X$ be a smooth Banach space and $T \in B(X)$. If $T$ is a surjective isometry, then $T$ is unitary. The converse holds when $X$ is reflexive.
\end{thm}
\begin{proof}
    Suppose that $T$ is a surjective isometry. Then for any nonzero $x \in X$,
    \[\|T'Jx\| \leq \|T'\|\|Jx\| = \|Jx\| = \|x\|.\]
    Now, using the surjectivity of $T$, choose a $y \in X$ such that $\displaystyle T(y)=\frac{x}{\|x\|}$. Note that $\|y\|=1$, $T$ being an isometry. This gives
    $$
    \|(T'Jx)(y)\| = \|Jx(Ty)\| = \frac{\|Jx(x)\|}{\|x\|}=\|x\|,
    $$
    proving that $ \|T'Jx\|\geq \|x\|$. Thus, we obtain $\|T^{\prime}Jx\|=\|x\|=\|Tx\|$.
    
    Conversely, if $T$ is unitary, then it is an isometry, by definition. Now, for any $f \in X^*$, there exists an $x\in X$ such that $Jx=f$, $J$ being surjective. So
    \[\|T'f\| = \|T'Jx\| = \|x\| = \|Jx\|=\|f\|.\]
    Hence $T'$ is bounded below, proving that $T$ is surjective (Theorem 13.9, \cite{limaye1996functional}).
\end{proof}

The next result, which establishes that the inverse of a unitary operator is also unitary, follows directly from Theorem \ref{unitary_charc}.

\begin{cor}\label{unitary_inv}
    Let $X$ be a smooth Banach space, and let $T \in B(X)$ be unitary and invertible. Then $T^{-1}$ is also a unitary operator on $X$.
\end{cor}

The following theorem provides a characterization of invertible unitary operators, and it also shows that the iso-abelian operators defined in \cite{stampfli1969}, coincide with the unitary operators defined here in the setting of smooth Banach spaces.

\begin{thm}\label{unitary_isoabelian}
    Let $X$ be a smooth Banach space and $T \in B(X)$ be invertible. Then  $T$ is unitary if and only if $J(Tx)=(T^{-1})'J(x)$ for each $x \in X$.
\end{thm}
\begin{proof}
    Let $T$ be a unitary operator on $X$. For $x \in X$, we have,
    $$(T^{-1})'J(x)(Tx)=J(x)(T^{-1}Tx)=\|x\|^2=\|Tx\|^2=J(Tx)(Tx).
    $$
    Also, since $T^{-1}$ is unitary, we get, $$\|(T^{-1})'J(x)\|=\|T^{-1}x\|=\|x\|=\|Tx\|.$$
    Since $X$ is smooth, we get $J(Tx)=(T^{-1})'J(x)$. For the converse part, for $x \in X$, we have, $$\|Tx\|^2=J(Tx)(Tx)=(T^{-1})'J(x)(Tx)=J(x)(x)=\|x\|^2.$$ Hence, $T$ is an isometry. Since $T$ is invertible, it is also onto. Therefore, by Theorem \ref{unitary_charc}, $T$ is a unitary operator.
\end{proof}

Using the above theorem, we can show that the composition of two invertible unitary operators is unitary.

\begin{prop}\label{unitary_composition}
    Let $T,S$ be invertible unitary operators on a smooth Banach space $X$. Then $T \circ S$ is a unitary operator. In particular, $T^n$ is unitary for every $n \in \mathbb{N}$.
\end{prop}
\begin{proof}
    For each $x\in X$, we have
    \[
(TS)^{-1\prime}J(x)=(T^{-1})'(S^{-1})'J(x)=(T^{-1})'J(Sx)=J(TSx).
\]
From Theorem \ref{unitary_isoabelian}, it follows that $T\circ S$ is a unitary operator.
\end{proof}

\begin{rem}
    For a self-adjoint operator $T\in B(X)$, it is known that $T^n$ is self-adjoint for each $n \in \mathbb{N}$, and that $T^{-1}$ is also self-adjoint whenever $T$ is invertible (Proposition 2.2, \cite{banachSA}). Analogous results hold for unitary operators (Propositions \ref{unitary_inv} and \ref{unitary_composition}). However, whether these properties hold for normal operators remains an open question.
\end{rem}
We now provide another characterization of unitary operators on reflexive smooth strictly convex Banach spaces, analogous to the Hilbert space case.

\begin{thm}
    Let $X$ be a reflexive smooth strictly convex Banach space. Then $T\in B(X)$ is unitary if and only if \[T^*T=TT^*=I.,\]
    where $T^*=J^{-1}T'J$.
\end{thm}
\begin{proof}
         Assume that $J^{-1}T'JT=TJ^{-1}T'J=I$. For any $x \in X$, we have
    \[\|Tx\|^{2} = J(Tx)(Tx) = T' J(Tx)(x) = J(x)(x) = \|x\|^{2}.\]
This proves that $T$ is an isometry. Now, since $TT^*=I$, $T$ is surjective, and hence, by Theorem \ref{unitary_charc}, $T$ is
a unitary operator.

    Conversely, if we assume that $T$ is an isometry, then for each $x \in X,$
    \[T' J(Tx)(x) = J(Tx)(Tx) = J(x)(x) = \|x\|^2.\]
    Now, for $x\neq 0$, we have, $\left\|T'J(Tx)\left(\frac{x}{\|x\|}\right)\right\| = \|x\|$, which implies $\|T'J(Tx)\| \geq \|x\|$. On the other hand,
    $$\|T'J(Tx)\| \leq \|T'\|\|JTx\| \leq \|JTx\|=\|Tx\|=\|x\|.$$
   Hence $T'J(Tx)$ satisfies $T' J(Tx)(x) = \|x\|^2$ and $\|T'J(Tx)\| = \|x\|$.  Using the uniqueness of the duality map, we get
   $$
    T'J(Tx) = J(x),
   $$
    for each $x\in X$, and hence, $J^{-1} T'JT = I$ (since $J$ is invertible). In addition, if $T$ is surjective, then for any $y \in X$ there exists an $x \in X$ such that $T(x) = y$. % Let $J(x) = x^{*}$ and $J(y) = y^{*}$. 
    Since $J^{-1} \, T' \, J \, T = I$, we get $T'(J(y)) = J(x),$ %that is, $T'(y^{*}) = x^{*}$,
    and therefore,
    \[
    T J^{-1} \, T' \, J(y) = T J^{-1}(Jx) = y.
    \]
    Hence, we obtain    \(T J^{-1} \, T' \, J = I,\) as required.
\end{proof}

We conclude this section by providing an example of a normal, self-adjoint and unitary operator on $l^p(\mathbb{N}), 1<p<\infty$.
\begin{ex}\label{examp}

Consider the operator $T$ on $\ell^p(\mathbb{N}), 1<p<\infty,$ defined by
$$T(x)=(x_2,x_1,x_3,x_4,\dots)\, \mbox{ for }\, x=(x_1,x_2,x_3,\dots) \in \ell^p(\mathbb{N}).$$
Consider the transpose $T'$ of this operator defined on $\ell^{q}(\mathbb{N})$ by
$$T'(y)=(y_2,y_1,y_3,\dots)\,\mbox{ for }\, y=(y_1,y_2,y_3,\dots)\in \ell^{q}(\mathbb{N}),$$
where $\frac1p+\frac1q=1$.
The normalized duality map $J:\ell^p(\mathbb{N}) \to \ell^{q}(\mathbb{N})$ is given as
$$J(x)=\frac{1}{\|x\|_p^{\,p-2}}
\big(|x_1|^{p-2}x_1,|x_2|^{p-2}x_2,|x_3|^{p-2}x_3,\dots\big).$$
Hence,
\[T'J(x) =\frac{1}{\|x\|_p^{\,p-2}}
\big(|x_2|^{p-2}x_2,|x_1|^{p-2}x_1,|x_3|^{p-2}x_3,\dots\big) 
\]
 and 
\[
J(Tx)
=
\frac{1}{\|x\|_p^{\,p-2}}
\big(|x_2|^{p-2}x_2,|x_1|^{p-2}x_1,\dots\big).
\]
Thus, $T'J(x)=J(Tx),$ so $T$ is self-adjoint, and hence a normal operator. On calculating the norm, we obtain
$$\|T'J(x)\|_{q}^{q}
=\frac{1}{\|x\|_p^{(p-2)q}}
\left(|x_2|^{(p-1)q}+|x_1|^{(p-1)q}+\cdots\right)
=\|x\|_p^{q}.$$
Since $T$ is an isometry on $\ell^p(\mathbb{N})$, we get
$$
\|T'J(x)\|_{q}=\|Tx\|_p =\|x\|_p \quad \text{for all } x\in \ell^p(\mathbb{N}),
$$
showing that $T$ is a unitary operator.

\end{ex}

In the following example shows that even if $T$ is normal, $T-\lambda I$ need not be normal.
\begin{ex}\label{normal_sum}
    Consider the same operator $T$ on $\ell^p(\mathbb{N}), 1<p<\infty, p\neq 2$ defined by
$$T(x)=(x_2,x_1,x_3,x_4,\dots)\, \mbox{ for }\, x=(x_1,x_2,x_3,\dots) \in \ell^p(\mathbb{N}).$$
Consider $e_1=(1,0,0,\ldots)$ and $\lambda =1$. Then $\|(T-I)e_1\| = \|(-1,1,0,0,\ldots)\|=2^{1/p}$. Since $J(e_1)=e_1$, similarly we get $\|(T-I)'Je_1\|=2^{1/q}$. These two quantities become same only if $p=2$. Hence $T-I$ is not a normal operator.
\end{ex}
Moreover, the above example shows that the sum of two normal, unitary, or self-adjoint operators need not be normal, unitary, or self-adjoint, respectively.

\section{Spectral properties of normal operators}
In this section, we establish some results on the spectrum of normal operators and unitary operators defined on smooth Banach spaces. In the Hilbert space setting, it is well known that for a self-adjoint operator, the eigenvalues are real, and the eigenvectors corresponding to distinct eigenvalues of normal operators are orthogonal. In the setting of smooth Banach spaces, it was shown (Theorem 6.3(i), \cite{banachSA}),  that the eigenvalues of self-adjoint operators are real. Here, we obtain a result analogous to the orthogonality of eigenvectors of normal operators in the setting of smooth Banach spaces. In order to establish this, we recall a known result from \cite{mal}.
\begin{prop}[Theorem 2.1.15, \cite{mal}]\label{birkoff}
    For any two elements $x,y$ in a Banach space $X$, $\|x+\lambda y\|\geq \|x\|$ for every $\lambda \in \C$ if and only if there exists a $\phi \in S_{X^*}$ such that $\phi(x)=\|x\|$ and $\phi(y)=0$.
\end{prop}

The following proposition on the connection between eigenvalues of a normal operator and its transpose helps to prove the orthogonality of the eigenvectors corresponding to distinct eigenvalues of normal operators on smooth Banach spaces and also 
establish the equality of the approximate spectrum for normal operators, in certain cases.

\begin{prop}\label{eigen_shift}
    Let $T \in B(X)$ be a normal operator and $\lambda \in \mathbb{C}$ be such that $Tx=\lambda x$ for some $x\neq 0$. Then $T'Jx=\lambda Jx$. That is, $\sigma_{eig}(T) \subseteq \sigma_{eig}(T')$. In addition, if $X$ is reflexive,
    \begin{enumerate}
        \item $\sigma_{eig}(T) = \sigma_{eig}(T')$
        \item $\sigma(T) = \sigma_{app}(T)$.
    \end{enumerate} 
\end{prop}
\begin{proof}
    Since $T$ is normal, $\|T'Jx\|=\|Tx\|=\|\overline{\lambda} x\|$, and
    $$T'J(x)(\overline{\lambda} x)=J(x)(\overline{\lambda} Tx)=\|\overline{\lambda} x\|^2.$$
    Since $X$ is smooth, this shows $T'Jx=\lambda Jx$. Hence $\sigma_{eig}(T) \subseteq \sigma_{eig}(T')$.

    Let $X$ be reflexive. Since $\sigma(T) = \sigma_{app}(T) \cup \sigma_{eig}(T')$, we need to prove $\sigma_{eig}(T')\subseteq \sigma_{eig}(T)$.  For $\lambda \in \sigma_{eig}(T')$, there exists a nonzero $f \in X^*$ such that $\|(T'-\lambda I)f\|=0$. Since $J$ is onto, we get a nonzero $x\in X$ such that $\|(T'-\lambda I)Jx\|=0$, giving $T'Jx = \lambda Jx$. Since $T'$ is normal (Proposition \ref{normal_transp}), using Proposition \ref{eigen_shift}, we get $T''\tilde{J}(Jx)=\lambda \tilde{J}(Jx)$, where $\tilde{J}$ is duality map of $X^*$. This implies $\lambda \in \sigma_{eig}(T'')=\sigma_{eig}(T)$
    
\end{proof}

As a direct consequence of the above result, we obtain the equality of the minimum modulus and Crawford number for a non-invertible normal operator, generalizing Proposition 3.15 of \cite{crawford2024}. 

\begin{cor}\label{non_inv_bound}
    Let $T$ be a normal operator on a reflexive smooth Banach space $X$. Then $T$ is invertible if and only if $T$ is bounded below. Hence, if $T$ is a non-invertible normal operator, then  $c(T)=\mu(T)=0.$
\end{cor}

Now we prove the orthogonality of eigenvectors corresponding to distinct eigenvalues of normal operators.

\begin{prop}\label{eigenvec_BKJ}
    Let $x_1$ and $x_2$ be eigenvectors corresponding to two distinct eigenvalues of a normal operator $T$ on $X$. Then $x_1$ and $x_2$ are Birkhoff James orthogonal to each other.
\end{prop}

\begin{proof}
    Let $Tx_1$ = $\lambda_1 x_1$ and $Tx_2=\lambda_2x_2$, where $\lambda_1 \neq \lambda_2$. Using Proposition \ref{eigen_shift}, we get, 
    $$\begin{aligned}
    \lambda_2J(x_1)(x_2) &= J(x_1)(\lambda_2x_2) =J(x_1)(Tx_2)\\
    &=T'J(x_1)(x_2)=\lambda_1J(x_1)(x_2).
    \end{aligned}
    $$
    Hence, if $\lambda_1-\lambda_2 \neq 0$, then $J(x_1)(x_2)=0$. Similarly, we get $J(x_2)(x_1)=0$. Hence, by Proposition \ref{birkoff}, we get the desired result.
\end{proof}

 It is known that the eigenspectrum of a normal operator on a separable Hilbert space is countable. The same holds for normal operators on separable smooth Banach spaces.

\begin{prop}\label{countable_SA}
    Let $X$ be a separable smooth Banach space, and $T \in B(X)$ be normal. Then $\sigma_{eig}(T)$ is countable.
\end{prop}
\begin{proof}
    Let $\Omega\subset X$ be the set constructed by taking one unit eigenvector corresponding to each distinct eigenvalue of $T$. Then for each $x\neq y$ from $\Omega$, by Proposition \ref{eigenvec_BKJ}, we have, $J(x)(y)=0$ and $\|J(x)\|=1=\|x\|$.  Now applying Proposition \ref{birkoff} with $\lambda=-1$, we get $\|x-y\|\geq 1$. Thus the elements of $\Omega$ are $1-$distance apart from each other.

    Now suppose that $\Omega$ is uncountable. Then consider the collection 
    $$A:=\left\{B(x, \frac{1}{2}):\, x\in \Omega\right\},$$ 
    where $B(a,r)$ represents an open ball centred at $a$ and of radius $r$. Then $A$ is an uncountable family of pairwise disjoint open balls in $X$. Consequently, any dense subset of the space must intersect each element of $A$, forcing every dense subset to be uncountable, which contradicts the separability of the space. Therefore, $\Omega$ must be countable.
\end{proof}

 In \cite{banachSA}, it was proved that for a self-adjoint operator $T$ on a smooth Banach space, $\|T^{2^n}\|=\|T\|^{2^n}$. In Theorem \ref{norm_spec_equal}, we establish a generalized version of this result for normal operators using the following theorem.

 \begin{thm}\label{normal_norm_ineq}
     Let $T \in B(X)$ be normal. Then $\|Tx\|^n\leq \|T^nx\|\|x\|^{n-1}$ for each $x \in X$ and $n \in \mathbb{N}$.
 \end{thm}
 \begin{proof}
     The result is trivial for $n=1$. Since $T$ is normal, we have $\|T^2x\|=\|T'J(Tx)\|$ for each $x\in S_X$. Moreover, $$\|Tx\|^2=J(Tx)(Tx)\leq \|T'J(Tx)\|\|x\|=\|T^2x\|\|x\|.$$ 
     Now, assume that $Tx\neq 0$ and $\|Tx\|^k\leq \|T^kx\|\|x\|^{k-1}$ for some $k\in \mathbb{N}$. We have
     $$
            \|Tx\|^{2k} = \|T^2x\|^{k}\|x\|^{k} \leq \|T^k(Tx)\|\|Tx\|^{k-1}\|x\|^k.
     $$
     Hence $\|Tx\|^{k+1} \leq \|T^{k+1}x\|\|x\|^k$, and the result follows from an induction argument.
 \end{proof}

The following corollary is immediate from the above theorem.

 \begin{cor}
     Let $T \in B(X)$ be a normal operator and $x \in X$. If $T^nx=0$ for some $n \in \mathbb{N}$, then $Tx=0$.
 \end{cor}

As a consequence of Theorem \ref{normal_norm_ineq}, we obtain relationships between norm, numerical radius, and minimum modulus for normal operators on smooth Banach spaces and their powers. Analogous results for self-adjoint operators on smooth Banach spaces were established in \cite{banachSA,sd2026}.

 \begin{thm}\label{norm_spec_equal}
    Let $T \in B(X)$ be a normal operator. Then,
    \begin{enumerate}
        \item $\|T^n\|=\|T\|^n$ for $n \in \mathbb{N}$.
        \item $\mu(T^{n})=\mu(T)^{n}$ for $n \in \mathbb{N}$.
        \item $r(T)=\rho(T)=\|T\|$.
        \item if $T \neq 0$ is compact, then there exists a $\lambda \in \sigma_{eig}(T)$ such that $|\lambda|=\|T\|$.
    \end{enumerate}
\end{thm}
\begin{proof}
    Since $T$ is normal, using Theorem \ref{normal_norm_ineq}, taking supremum over $x \in S_X$, we get $\|T^{n}\|=\|T\|^n$ for each $n \in \mathbb{N}$. A similar argument works for the minimum modulus. Using the spectral radius formula, we get $\rho(T)=\|T\|$. Since $\sigma(T) \subseteq \overline{W(T)}$, we get $r(T)=\rho(T)=\|T\|$. The spectrum being compact, there exist a $\lambda \in \sigma(T)$ such that $|\lambda|=\rho(T)=\|T\|$. For a nonzero compact operator $T$, since every non-zero spectral value is an eigenvalue, we obtain an eigenvalue $\lambda$ such that $|\lambda|=\|T\|$.
\end{proof}

In the Hilbert space setting, for a normal operator $T$ and $k\notin\sigma(T)$, it is well known that $\displaystyle\|(T-kI)^{-1}\|=\frac{1}{\operatorname{dist}(k,\sigma(T))}$ (\cite{limaye1996functional}). A partial extension of this result is established on smooth Banach spaces.

\begin{prop}\label{norm_spec_dist}
    Let $T \in B(X)$ be such that $(T-kI)^{-1}$ is a  normal operator for some $k \notin \sigma(T)$. Then,
    $$
    \|(T-kI)^{-1}\|=\frac{1}{dist(k,\sigma(T))}.
    $$
    In particular, if $T \in B(X)$ is invertible and $T^{-1}$ is normal, then,
    $$
    \frac{1}{\|T^{-1}\|} = \min_{\lambda \in \sigma(T)}|\lambda|. 
    $$
    
\end{prop}

\begin{proof}
    Since $(T-kI)^{-1}$ is normal, by Theorem \ref{norm_spec_equal}, we have $\|(T-kI)^{-1}\|=\rho((T-kI)^{-1})$. Now, appying Theorem 10.2.12, \cite{MTNfunctional}, we get
    $$
    \|(T-kI)^{-1}\|=\frac{1}{dist(k,\sigma(T))}.
    $$
    Since $\sigma(T)$ is compact, the particular case follows by choosing $k=0$.
\end{proof}

For a subset $E$ of a smooth Banach space $X$, we use the notation $E^\perp$ to denote the set $\{y\in X: J(x)(y) =0\, \forall \, x \in E \}$. Then we have following proposition.

\begin{prop}
    Let $T$ be a normal operator on $X$. Let $E$ be set of all eigenvectors of $T$. Then $\overline{spanE}$ and $E^\perp$ are closed invariant subspaces of $X$ under $T$.
\end{prop}

\begin{proof}
    Clearly, $spanE$ is an invariant subspace of $T$. Therefore $\overline{spanE}$ is a closed invariant subspace of $X$. Let $y \in E^\perp.$ That is $J(x)(y)=0$ for all $x \in E$. Now, for each $x \in E$, we have $Tx = \lambda x$. Then using proposition \ref{eigen_shift}, we have
    $$
    J(x)(Ty)=T'J(x)(y) = \lambda J(x)(y)=0.
    $$
    Hence $Ty \in E^\perp$. Since $ E^\perp=\displaystyle\bigcap_{x\in E} kerJ(x)$, $E^\perp$ is a closed invariant subspace for $T$.
\end{proof}

\begin{rem}
    Since every self-adjoint and unitary operator on $X$ are normal, the results established for normal operators in this and the previous section apply to self-adjoint and unitary operators also. Some of these results have already been established for self-adjoint operators in \cite{banachSA,sd2026,stampfli1969}.

\end{rem}

The spectrum of a unitary operator on a smooth Banach space is quite analogous to that of its Hilbert space counterpart.
\begin{prop}\label{unit_spec}
    Let $T \in B(X)$ be an invertible unitary operator. Then $\sigma(T)$ contained in the unit circle of $\C$.
\end{prop}
\begin{proof}
    Since $T$ and $T^{-1}$ are isometries $\|T\|=\|T^{-1}\|=1$. This implies that $\sigma(T)$ is contained in the closed unit ball centered at origin. Now, for $\lambda$ with $|\lambda|<1$, we have $||\lambda T^{-1}\|<1$. Hence $I-\lambda T^{-1}$ is invertible. But $T-\lambda I = T(I-\lambda T^{-1})$, showing that $\lambda \notin \sigma(T)$ for  $|\lambda|<1$. Hence for $\lambda \in \sigma(T)$, we have $|\lambda|=1$.
\end{proof}

\begin{rem}
    When $X$ is a reflexive smooth Banach space, every unitary operator is an onto isometry and hence is invertible (Theorem \ref{unitary_charc}). Consequently, in the setting of reflexive smooth Banach spaces, the assumption of invertibility can be omitted in Proposition \ref{unitary_inv}, Theorem \ref{unitary_isoabelian}, Proposition \ref{unitary_composition} and Proposition \ref{unit_spec}.
\end{rem}

\section*{Conclusion and open problems}

In this article, we have introduced classes of normal and unitary operators on smooth Banach spaces and investigated several of their spectral properties. In particular, we obtained a number of results that are analogous to their counterparts in the Hilbert space setting. At the same time, the extent to which other well-known properties from Hilbert space theory remain valid in the setting of smooth Banach spaces is not yet fully clear and needs further investigation. The following are some of the natural questions arise in this regard. 
\begin{enumerate}
    
\item What can be said about the structure of the numerical range of special operators such as normal, unitary, and self-adjoint operators? More specifically, under what conditions  the numerical range is convex? 
\item Let $X$ be a reflexive strictly convex smooth Banach space and $T \in B(X)$. Is it true that $T$ is normal if and only if $TT^*=T^*T$, where $T^*=J^{-1}T'J$.
\end{enumerate}

 \section*{Acknowledgments}
 The authors would like to thank Dr. C. R. Jayanarayanan, Associate Professor, Indian Institute of Technology, Palakkad, for his valuable suggestions and insightful comments. Mr. Mohammed Shameem would like to thank the University Grants Commission of India for funding his research work (Fellowship No. 231610086159).

\section*{Statements and Declarations}
 The corresponding author declares that there is no conflict of interest on behalf of all authors. During this study, no datasets were created or examined.

\end{document}